\documentclass[12pt,oneside,reqno]{amsart}

\usepackage{amssymb}
\usepackage{algpseudocode}
\usepackage{algorithm}
\usepackage{verbatim}
\usepackage{graphicx}
\usepackage{placeins}
\usepackage{listings}
\usepackage{xcolor}
\usepackage{subcaption}
\usepackage[noadjust]{cite}
\RequirePackage{dsfont} 

 \scrollmode
\allowdisplaybreaks
\usepackage{graphicx}
\usepackage{epstopdf}
\usepackage{hyperref}

\usepackage{amsmath}
\usepackage{tikz}

\makeatletter
\newcommand{\xleftrightarrow}[2][]{\ext@arrow 3359\leftrightarrowfill@{#1}{#2}}
\makeatother

\definecolor{codegreen}{rgb}{0,0.6,0}
\definecolor{codegray}{rgb}{0.5,0.5,0.5}
\definecolor{codepurple}{rgb}{0.58,0,0.82}
\definecolor{backcolour}{rgb}{0.95,0.95,0.92}

\lstdefinestyle{list_style}{
  backgroundcolor=\color{backcolour}, commentstyle=\color{codegreen},
  keywordstyle=\color{magenta},
  numberstyle=\tiny\color{codegray},
  stringstyle=\color{codepurple},
  basicstyle=\ttfamily\footnotesize,
  breakatwhitespace=false,         
  breaklines=true,                 
  captionpos=b,                    
  keepspaces=true,                 
  numbers=left,                    
  numbersep=5pt,                  
  showspaces=false,                
  showstringspaces=false,
  showtabs=false,                  
  tabsize=2
}

\newcommand{\xdasharrow}[2][->]{
\tikz[baseline=-\the\dimexpr\fontdimen22\textfont2\relax]{
\node[anchor=south,font=\scriptsize, inner ysep=1.5pt,outer xsep=2.2pt](x){#2};
\draw[shorten <=3.4pt,shorten >=3.4pt,dashed,#1](x.south west)--(x.south east);
}
}

\newcommand{\DEBUG}{}

\ifdefined\DEBUG%

  \def\rem#1{{\marginpar{\raggedright\scriptsize #1}}}
  \newcommand{\pmr}[1]{\rem{\color{blue}{$\bullet$ #1}}}
  \newcommand{\ppr}[1]{\rem{\color{red}{$\bullet$ #1}}}
 \else%

  \newcommand{\ppr}[1]{}
  \newcommand{\pmr}[1]{}
 \fi

\newtheorem{theorem}{Theorem}[section]

\newenvironment{pf*}[1]{\noindent\textit{#1} }{}

\begin{document}

\title
[Monte Carlo with adaptive variance reduction]
{Monte Carlo integration of $C^r$ functions with adaptive variance reduction: an asymptotic analysis}

\author[L. Plaskota]{Leszek Plaskota}
\address{University of Warsaw, 
Faculty of Mathematics, Informatics and Mechanics, ul. S.~Banacha 2, 02-097 Warsaw, Poland}
\email{leszekp@mimuw.edu.pl, corresponding author}

\author[P. Przyby{\l}owicz]{Pawe{\l} Przyby{\l}owicz}
\address{AGH University of Science and Technology,
Faculty of Applied Mathematics,
Al. A.~Mickiewicza 30, 30-059 Krak\'ow, Poland}
\email{pprzybyl@agh.edu.pl}

\author[{\L}. Stepie\'{n}]{{\L}ukasz St\c epie\'{n}}
\address{AGH University of Science and Technology,
	Faculty of Applied Mathematics,
	Al. A.~Mickiewicza 30, 30-059 Krak\'ow, Poland}
\email{lstepie@agh.edu.pl}

\begin{abstract}
The theme of the present paper is numerical integration of $C^r$ functions using randomized methods. 
We consider variance reduction methods that consist in two steps. First the initial interval is partitioned into 
subintervals and the integrand is approximated by a piecewise polynomial interpolant that is based on 
the obtained partition. Then a randomized approximation is applied on the difference of the integrand and 
its interpolant. The final approximation of the integral is the sum of both. The optimal convergence rate is 
already achieved by uniform (nonadaptive) partition plus the crude Monte Carlo; however, special adaptive 
techniques can substantially lower the asymptotic factor depending on the integrand. The improvement 
can be huge in comparison to the nonadaptive method, especially for functions with rapidly varying 
$r$th derivatives, which has serious implications for practical computations. In addition, the proposed 
adaptive methods are easily implementable and can be well used for automatic integration.
\end{abstract}

\maketitle


\section{Introduction}

Adaption is a useful tool to improve performance of algorithms. The problems of 
numerical integration and related to it $L^1$ approximation are not exceptions, 
see, e.g.,~\cite{Novak1996} for a survey of theoretical results on the subject.
If an underlying function possesses some singularities and is otherwise smooth, 
then using adaption is necessary to localize the singular points and restore 
the convergence rate typical for smooth functions, see, e.g., 
\cite{KP1, PlaskotaWasilkowski2005, PlaskotaWasilkowskiZhao2008, PP1}. 
For functions that are smooth in the whole domain, adaptive algorithms usually do not offer 
a better convergence rate than nonadaptive algorithms; however, they can essentially 
lower asymptotic constants. This is why adaptive quadratures are widely used for numerical 
integration, see, e.g., \cite{Gonnet2012,Lyness1972, DavisRabinowitz1984}. Their superiority 
over nonadaptive quadratures is rather obvious, but precise answers to the quantitative 
question of "how much adaption helps" are usually missing. This gap was partially filled by 
recent results of 
\cite{Gocwin2021, Plaskota2015, PlaskotaSamoraj2022}, 
where best asymptotic constants of deterministic algorithms that use piecewise 
polynomial interpolation were determined for integration and $L^1$ approximation of 
$r$-times continuously differentiable functions $f:[a,b]\to\mathbb R.$ 
In this case, adaption relies on adjusting 
the partition of the interval $[a,b]$ to the underlying function. 
While the convergence rate is of order $N^{-r},$ where $N$ is the number of 
function evaluations used, it turns out that 
the asymptotic constant depends on $f$ via the factor of 
$(b-a)^r\big\|f^{(r)}\big\|_{L^1}$ for uniform partition, 
and $\big\|f^{(r)}\|_{L^{1/(r+1)}}$ for best adaptive partition. 

In the current paper, we present the line of thinking similar to that of 
the aforementioned papers. The difference is that now we want to carry out
the analysis and obtain asymptotic constants for randomized algorithms.

Our goal is the numerical approximation of the integral 
\begin{equation}\label{theproblem} 
Sf=\int_a^bf(x)\,\mathrm dx.
\end{equation}
It is well known that for $f\in L^2(a,b)$ the \emph{crude Monte Carlo},
\begin{equation}\label{MCstandard}
 M_Nf=\frac{b-a}N\sum_{i=1}^N f(t_i),\quad\mbox{where}\quad t_i\stackrel{iid}\sim U(a,b),
 \end{equation} 
returns an approximation with expectation $\mathbb E(M_Nf)=Sf$ and error (standard deviation)
\begin{equation}\label{MCerr}
 \sqrt{\mathbb E\big(Sf-M_Nf\big)^2}=\frac{\sigma(f)}{\sqrt N},\quad\mbox{where}\quad
\sigma(f)^2=(b-a)S(f^2)-(Sf)^2.
\end{equation}
Suppose that the function enjoys more smoothness, i.e., $$f\in C^r([a,b]).$$ 
Then a much higher convergence rate $N^{-(r+1/2)}$ can be achieved using various 
techniques of \emph{variance reduction}, see, e.g., \cite{Heinrich1993}. One way is to apply 
a randomized approximation of the form
\begin{equation}\label{MCvr}
\overline M_{N,r}(f)=S(L_{m,r}f)+M_n(f-L_{m,r}f),
\end{equation} 
where $L_{m,r}$ is the piecewise polynomial interpolation of $f$ of degree $r-1$ using 
a partition of the interval $[a,b]$ into $m$ subintervals, $M_n$ is a Monte Carlo type 
algorithm using $n$ random samples of $f,$ and $N$ is the total number of function evaluations 
used (for arguments chosen either deterministically or randomly). 
The optimal rate is already achieved for uniform (nonadaptive) partition and 
crude Monte Carlo. Then, see Theorem~\ref{thm:equ} with $\beta=0,$ 
the error asymptotically equals 
$$c\,(b-a)^{r+1/2}\big\|f^{(r)}\big\|_{L^2(a,b)}\,N^{-(r+1/2)},$$
where $c$ depends only on the choice of the interpolation points within subintervals. 
The main result of this paper relies on showing that with the help of adaption 
the asymptotic error of the methods \eqref{MCvr} can be reduced to
\begin{equation}\label{bestada}
c\,\big\|f^{(r)}\big\|_{L^{1/(r+1)}(a,b)}\,N^{-(r+1/2)},
\end{equation}
see Theorems \ref{thm:strata} and \ref{thm:importstrata}. 
Observe that the gain can be significant, especially when the derivative $f^{(r)}$ 
drastically changes. For instance, for $r=4,$ $[a,b]=[0,1],$ and $f(x)=1/(x+d),$ 
adaption is asymptotically roughly $5.7*10^{12}$ times better than nonadaption 
if $d=10^{-4},$ and $1.8*10^{29}$ times if $d=10^{-8}.$ 

We construct two randomized algorithms, denoted $\overline M_{N,r}^{\,*}$ and 
$\overline M_{N,r}^{\,**},$ that achieve the error \eqref{bestada}. Although they  
use different initial approaches; namely, stratification versus importance sampling, 
in the limit they both reach essentially the same partition, such that the $L^1$ errors 
of Lagrange interpolation in all subintervals are equalized. However, numerical tests 
of Section \ref{sec:impl} show that the algorithm $\overline M_{N,r}^{\,*}$ achieves 
the error \eqref{bestada} with some delay, which makes $\overline M_{N,r}^{\,**}$ 
worth recommending rather than $\overline M_{N,r}^{\,*}$ in practical computations.

Other advantages of $\overline M_{N,r}^{\,**}$ are that it is easily implementable and,
as shown in Section \ref{sec:automatic}, it can be successfully used for automatic 
Monte Carlo integration.

Our analysis has been so far restricted to one-dimensional integrals only. In a~future 
work it will be extended and corresponding adaptive Monte Carlo algorithms will 
be constructed for multivariate integration, where randomization finds its major 
application. The current paper is the first step in this direction.

In the sequel, we use the following notation. For two functions of $N$ we write 
$g_1(N)\lessapprox g_2(N)$ iff $\limsup_{N\to\infty}g_1(N)/g_2(N)\le 1,$ and we write 
$g_1(N)\approx g_2(N)$ iff $\lim_{N\to\infty}g_1(N)/g_2(N)=1.$ Similarly, for functions 
of $\varepsilon$ we write $h_1(\varepsilon)\lessapprox h_2(\varepsilon)$ iff 
$\limsup_{\varepsilon\to 0^+}h_1(\varepsilon)/h_2(\varepsilon)\le 1,$ and 
$h_1(\varepsilon)\approx h_2(\varepsilon)$ iff 
$\lim_{\varepsilon\to 0^+}h_1(\varepsilon)/h_2(\varepsilon)=1.$

\section{Variance reduction using Lagrange interpolation}\label{sec:equispaced}

We first derive some general error estimates for the variance reduction algorithms of 
the form \eqref{MCvr}, where the standard Monte Carlo is applied for the error of 
piecewise Lagrange interpolation. Specifically, we divide the interval $[a,b]$ into $m$ 
subintervals using a partition $a=x_0<x_1<\cdots<x_m=b,$ and on each subinterval 
$[x_{j-1},x_j]$ we approximate $f$ using Lagrange interpolation of degree $r-1$ with 
the interpolation points 
\begin{equation*}
x_{j,s}=x_{j-1}+z_s(x_j-x_{j-1}),\qquad 1\le s\le r,
\end{equation*} where
\begin{equation}\label{zi-points}
0\le z_1<z_2<\cdots<z_r\le 1
\end{equation}
are fixed (independent of the partition). Denote such an approximation by $L_{m,r}f.$ 
Then $f=L_{m,r}f+R_{m,r}f$ with $R_{m,r}f=f-L_{m,r}f.$ 
The integral $Sf$ is finally approximated by
$$ \overline M_{m,n,r}f\,=\,S(L_{m,r}f)+M_n(R_{m,r}f), $$ 
where $M_n$ is the standard Monte Carlo \eqref{MCstandard}. We obviously have 
$\mathbb{E}(\overline M_{m,n,r}f)=Sf.$ Since 
$$ Sf-\overline M_{m,n,r}f\,=\,Sf-S(L_{m,r}f)-M_n(R_{m,r}f)\,=\,S(R_{m,r}f)-M_n(R_{m,r}f),$$
by \eqref{MCerr} we have
\begin{equation*}
\mathbb E\big(Sf-\overline M_{m,n,r}f\big)^2\,=\,
     \frac1n\left((b-a)S\big((R_{m,r}f)^2\big)-\big(S(R_{m,r}f)\big)^2\right).
\end{equation*}

\smallskip
Note that
$$ S\big((R_{m,r}f)^2\big)\,=\,\int_a^b(f-L_{m,r}f)^2(x)\,\mathrm dx
    \,=\,\|f-L_{m,r}f\|_{L^2(a,b)}^2 $$
is the squared $L^2$-error of the applied (piecewise) polynomial interpolation, while 
$$ S(R_{m,r}f)\,=\,\int_a^b(f-L_{m,r}f)(x)\,\mathrm dx\,=\,S(f)-S(L_{m,r}f) $$
is the error of the quadrature $\overline Q_{m,r}f=S(L_{m,r}f).$ 

\smallskip
From now on we assume that $f$ is not a polynomial of degree smaller than or equal to 
$r-1,$ since otherwise $\overline M_{m,n,r}f=Sf.$ Define the polynomial
\begin{equation}\label{mainpoly}
P(z)=(z-z_1)(z-z_2)\cdots(z-z_r).
\end{equation}

\medskip
We first consider the interpolation error $\|f-L_{m,r}f\|_{L^2(a,b)}.$ Let
\begin{equation}\label{alpha}
  \alpha\,=\,\|P\|_{L^2(0,1)}=\bigg(\int_0^1|P(z)|^2\mathrm dz\bigg)^{1/2}. 
\end{equation}
For each $j,$ the local interpolation error equals
\begin{eqnarray*} \big\|f-L_{m,r}f\big\|_{L^2(x_{j-1},x_j)} &=&
   \bigg(\int_{x_{j-1}}^{x_j}\bigl|\,(x-x_{j,1})\cdots(x-x_{j,r})f[x_{j,1},\ldots,x_{j,r},x]\,\bigr|^2\mathrm dx\bigg)^{1/2} \\
   &=& \alpha\,h_j^{r+1/2}\,\frac{|f^{(r)}(\xi_j)|}{r!},\qquad\qquad\xi_j\in[x_{j-1},x_j].
\end{eqnarray*}
Hence 
$$ \|f-L_{m,r}f\|_{L^2(a,b)}\,=\,\frac{\alpha}{r!}\bigg(\sum_{j=1}^m h_j^{2r+1}\big|f^{(r)}(\xi_j)\big|^2\bigg)^{1/2}.$$
In particular, for the equispaced partition, in which case $h_j=(b-a)/m,$ we have
\begin{eqnarray*} 
     \big\|f-L_{m,r}f\big\|_{L^2(a,b)} &=& \frac{\alpha}{r!}\,\bigg(\frac{b-a}{m}\bigg)^r
     \bigg(\frac{b-a}{m}\sum_{j=1}^m|f^{(r)}(\xi_j)|^2\bigg)^{1/2} \\
   &\approx& \frac{\alpha}{r!}\,\bigg(\frac{b-a}{m}\bigg)^r\,\big\|f^{(r)}\big\|_{L^2(a,b)}
     \qquad\mbox{as}\quad m\to+\infty. \nonumber
\end{eqnarray*}

\medskip
Now, we consider the quadrature error $Sf-\overline Q_{m,r}f.$ Let 
\begin{equation}\label{beta}
 \beta\,=\,\int_0^1 P(z)\,\mathrm dz. 
\end{equation}
The local integration errors equal
\begin{eqnarray*}
  \lefteqn{\int_{x_{j-1}}^{x_j}(f-L_{m,r}f)(x)\,\mathrm dx \;=\;
   \int_{x_{j-1}}^{x_j} (x-x_{j,1})\cdots(x-x_{j,r})f[x_{j,1},\ldots,x_{j,r},x]\,\mathrm dx}  \\
   &&=\;\frac 1{r!}\int_{x_{j-1}}^{x_j} (x-x_{j,1})\cdots(x-x_{j,r})f^{(r)}(\xi_j(x))\,\mathrm dx,
    \qquad\xi_j(x)\in[x_{j-1},x_j].
\end{eqnarray*} 
Choose arbitrary $\zeta_j\in[x_{j-1},x_j]$ for $1\le j\le m.$ Then 
\begin{eqnarray*} 
  \lefteqn{\bigg|\frac 1{r!}\,\int_{x_{j-1}}^{x_j} (x-x_{j,1})\cdots(x-x_{j,r})f^{(r)}(\xi_j(x))\,\mathrm dx\,-\,
     \frac{f^{(r)}(\zeta_j)}{r!}\int_{x_{j-1}}^{x_j} (x-x_{j,1})\cdots(x-x_{j,r})\,\mathrm dx\bigg| } \\ 
     &=& \frac1{r!}\,\bigg|\int_{x_{j-1}}^{x_j}(x-x_{j,1})\cdots(x-x_{j,r})\left(f^{(r)}(\xi_j(x))-f^{(r)}(\zeta_j)\right)\,
      \mathrm dx\bigg|\;\le\;\omega(h_j)\,\frac{h_j^{r+1}}{r!}\,\|P\|_{L^1(0,1)},
\end{eqnarray*}
where $\omega$ is the modulus of continuity of $f^{(r)}.$ We also have 
$$ \frac{f^{(r)}(\zeta_j)}{r!}\int_{x_{j-1}}^{x_j}(x-x_{j,1})\cdots(x-x_{j,r})\,\mathrm dx\,=\,
      \frac{\beta}{r!}\,h_j^{r+1}f^{(r)}(\zeta_j). $$
Hence $Sf-\overline Q_{m,r}f\,=\,X_m\,+\,Y_m,$ where
$$  X_m \,=\, \frac{\beta}{r!}\,\sum_{j=1}^mh_j^{r+1}f^{(r)}(\zeta_j)\qquad\mbox{and}\qquad
    |Y_m| \,\le\, \frac{\|P\|_{L^1(0,1)}}{r!}\sum_{j=1}^m\omega(h_j)h_j^{r+1}. $$
In particular, for the equispaced partition,
\begin{eqnarray*}
    X_m &=& \frac{\beta}{r!}\,(b-a)^r\bigg(\sum_{j=1}^m\frac{b-a}{m}f^{(r)}(\zeta_j)\bigg)\,m^{-r}, \\
  |Y_m| &\le& \frac{\|P\|_{L^1(0,1)}}{r!}\;\omega\bigg(\frac{b-a}{m}\bigg)(b-a)^{r+1}m^{-r}.
\end{eqnarray*}

Suppose that $\beta\ne 0$ and $\int_a^bf^{(r)}(x)\,\mathrm dx\ne 0.$ Then
$X_m\approx\frac{\beta}{r!}(b-a)^r\left(\int_a^bf^{(r)}(x)\,\mathrm dx\right)m^{-r}.$
Since $\omega(h)$ goes to zero as $h\to 0^+,$ the component $X_m$ dominates $Y_m$ 
as $m\to+\infty.$ Hence
\begin{equation*}
   Sf-\overline Q_{m,r}f \,\approx\, 
   \frac{\beta}{r!}\,\bigg(\frac{b-a}{m}\bigg)^r\,\bigg(\int_a^bf^{(r)}(x)\,\mathrm dx\bigg)
     \qquad\mbox{as}\quad m\to+\infty.
\end{equation*}
On the other hand, if $\beta=0$ or $\int_a^bf^{(r)}(x)\,\mathrm dx=0$ then the quadrature 
error converges to zero faster than $m^{-r},$ i.e.
$$\lim_{m\to+\infty}\big(Sf-\overline Q_{m,r}f\big)\,m^r\,=\,0.$$ 
Note that $\beta=0$ if and only if the quadrature $\overline Q_{m,r}$ has the degree of exactness 
at least $r,$ i.e., it is exact for all polynomials of degree $r$ or less. Obviously, the maximal degree 
of exactness equals $2r-1.$ 

\medskip
We see that for the equidistant partition of the interval $[a,b]$ the error 
$\big(\mathbb E(Sf-\overline M_{m,n,r}f)^2\big)^{1/2}$ is asymptotically proportional 
to $$\phi(m,n)=n^{-1/2}m^{-r},$$  
regardless of the choice of points $z_i$ in \eqref{zi-points}. 
Let us minimize $\phi(m,n)$ assuming the total number of points used is at most $N.$ 
We have two cases depending on whether both endpoints of each subinterval are used 
in interpolation. If so, i.e., if $z_1=0$ and $z_r=1$ (in this case $r\ge 2$) then 
$N=(r-1)m+1+n.$ The optimal values are 
\begin{equation}\label{mn1}
m^*=\frac{2r(N-1)}{(r-1)(2r+1)},\qquad n^*=\frac{N-1}{2r+1},
\end{equation} 
for which $$\phi(m^*,n^*)\,=\,
\sqrt{2}\,\bigg(1-\frac1r\bigg)^r\bigg(\frac{r+1/2}{N}\bigg)^{r+1/2}.$$
Otherwise we have $N=rm+n.$ The optimal values are 
\begin{equation}\label{mn2}
m^*=\frac{2N}{2r+1},\qquad n^*=\frac{N}{2r+1},
\end{equation} for which
$$\phi(m^*,n^*)\,=\,\sqrt{2}\,\bigg(\frac{r+1/2}{N}\bigg)^{r+1/2}.$$

\smallskip 
Denote by $\overline M_{N,r}$ the corresponding algorithm with
the equidistant partition, where for given $N$ the values 
of $n$ and $m$ equal correspondingly $\lfloor n^*\rfloor$ and $\lfloor m^*\rfloor.$ 
Our analysis is summarized in the following theorem.
\begin{theorem}\label{thm:equ} We have as $N\to+\infty$ that
$$ \sqrt{\mathbb E\big(Sf-\overline M_{N,r}f\big)^2}\;\approx\;
c_r\,(b-a)^r\,C(P,f)\,N^{-(r+1/2)},$$ where
$$ C(P,f)=\sqrt{\alpha^2\,(b-a)\bigg(\int_a^b\big|f^{(r)}(x)\big|^2\mathrm dx\bigg)\,-\,
      \beta^2\bigg(\int_a^bf^{(r)}(x)\,\mathrm dx\bigg)^2},$$ 
$\alpha$ and $\beta$ are given by \eqref{alpha} and \eqref{beta}, and
\begin{equation}\label{ciar}
c_r=\left\{\begin{array}{ll}\sqrt{2}\,\big(1-\frac1r\big)^r\frac{(r+1/2)^{r+1/2}}{r!},&
\quad\mbox{if}\quad r\ge 2,\,z_1=0,\,z_r=1,\\ \ \sqrt{2}\,\frac{(r+1/2)^{r+1/2}}{r!},&
\quad\mbox{otherwise}.\end{array}\right.
\end{equation}
\end{theorem}

We add that the algorithm $\overline M_{N,r}$ is fully implementable since we assume that 
we have access to function evaluations at points from $[a,b].$

\section{First adaptive algorithm}\label{sec:strata}

Now we add a stratification strategy to our algorithm of 
Theorem \ref{thm:equ} to obtain an adaptive algorithm with 
a much better asymptotic constant. That is, 
we divide the initial interval $[a,b]$ into $k$ equal length subintervals $I_i, \ 1\leq i \leq k,$ and on each 
subinterval we apply the approximation of Theorem \ref{thm:equ} with some $N_i,$ where 
\begin{equation}\label{Nsum}\sum_{i=1}^k N_i\le N.\end{equation}
Denote such an approximation by $\overline M_{N,k,r}.$ 
(Note that $\overline M_{N,r}=\overline M_{N,1,r}$.) Then, by Theorem \ref{thm:equ}, 
for fixed $k$ we have as all $N_i\to+\infty$ that
$$ \sqrt{\mathbb E\big(Sf-\overline M_{N,k,r}f\big)^2}\;\approx\;c_r h^r
     \bigg(\sum_{i=1}^k\frac{C_i^2}{N_i^{2r+1}}\bigg)^{1/2},$$ where  
\begin{equation*}
C_i=C_i(P,f)=\sqrt{\alpha^2\,h\,\int_{I_i}\big|f^{(r)}(x)\big|^2\mathrm dx-\beta^2\,
\left(\int_{I_i}f^{(r)}(x)\,\mathrm dx\right)^2},\qquad h=\frac{b-a}{k}.
\end{equation*}
Minimization of $\psi(N_1,\ldots,N_k)=\left(\sum_{i=1}^kC_i^2N_i^{-(2r+1)}\right)^{1/2}$ with respect 
to \eqref{Nsum} gives
$$ N_i^*\,=\,\frac{C_i^{1/(r+1)}}{\sum_{j=1}^k C_j^{1/(r+1)}}\,N,\qquad 1\le i\le k, $$
and then 
$$\psi(N_1^*,\ldots,N_k^*)= \bigg(\sum_{i=1}^k C_i^{1/(r+1)}\bigg)^{r+1}N^{-(r+1/2)}. $$

Let $\xi_i,\eta_i\in I_i$ satisfy 
$\int_{I_i}\big|f^{(r)}(x)\big|^2\mathrm dx=h\big|f^{(r)}(\xi_i)\big|^2$ and 
$\int_{I_i}f^{(r)}(x)\,\mathrm dx=hf^{(r)}(\eta_i).$ Then 
$$C_i=h\sqrt{\alpha^2|f^{(r)}(\xi_i)|^2-\beta^2|f^{(r)}(\eta_i)|^2}$$ 
and we have as $k\to+\infty$ that
\begin{eqnarray*}
    \bigg(\sum_{i=1}^kC_i^{1/(r+1)}\bigg)^{r+1} &=& h\,\bigg(\sum_{i=1}^k
    \big(\alpha^2|f^{(r)}(\xi_i)|^2-\beta^2|f^{(r)}(\eta_i)|^2\big)^\frac{1}{2(r+1)}\bigg)^{r+1}\nonumber \\
    &\approx&h\,(\alpha^2-\beta^2)^{1/2}\bigg(\sum_{i=1}^k\big|f^{(r)}(\xi_i)\big|^{1/(r+1)}\bigg)^{r+1}\nonumber \\
    &\approx&h^{-r}(\alpha^2-\beta^2)^{1/2}\bigg(\sum_{i=1}^kh\big|f^{(r)}(\xi_i)\big|^{1/(r+1)}\bigg)^{r+1}\nonumber \\
    &\approx&h^{-r}(\alpha^2-\beta^2)^{1/2}\bigg(\int_a^b\big|f^{(r)}(x)\big|^{1/(r+1)}\mathrm dx\bigg)^{r+1}.
\end{eqnarray*}

It is clear that we have to take $N_i$ to be an integer and at least $r,$ for instance 
$$ N_i=\left\lfloor N_i^*\left(1-\frac{kr}N\right)+r\right\rfloor,\qquad 1\le i\le k. $$
Then the corresponding number $m_i$ of subintervals and number $n_i$ of random points in $I_i$
can be chosen as
$$ m_i=\max\left(\lfloor m_i^*\rfloor,1\right),\qquad n_i=\lfloor n_i^*\rfloor,$$
where $m_i^*$ and $n_i^*$ are given by \eqref{mn1} and \eqref{mn2} with $N$ replaced by $N_i.$

Denote by $\overline M^{\,*}_{N,r}$ the above constructed approximation $\overline M_{N,k_N,r}$ with 
$k_N$ such that $k_N\to+\infty$ and $k_N/N\to 0$ as $N\to+\infty.$ For instance, $k_N=N^\kappa$ with 
$0<\kappa<1.$ Our analysis gives the following result.

\begin{theorem}\label{thm:strata} 
We have as $N\to+\infty$ that
$$ \sqrt{\mathbb E\big(Sf-\overline M^{\,*}_{N,r}f\big)^2}\,\approx\,
    c_r\,\sqrt{\alpha^2-\beta^2}\,
    \bigg(\int_a^b\big|f^{(r)}(x)\big|^{1/(r+1)}\mathrm dx\bigg)^{r+1}N^{-(r+1/2)}.$$
\end{theorem}

The asymptotic constant of the approximation $\overline M^{\,*}_{N,r}$ of Theorem \ref{thm:strata} 
is never worse than that of $\overline M_{N,r}$ of Theorem \ref{thm:equ}.
Indeed, comparing both constants we have
\begin{eqnarray*}
&&c_r(b-a)^r\sqrt{\alpha^2\,(b-a)\bigg(\int_a^b\big|f^{(r)}(x)\big|^2\mathrm dx\bigg)\,-\,
      \beta^2\bigg(\int_a^bf^{(r)}(x)\,\mathrm dx\bigg)^2}\\
      &&\qquad\ge\;c_r\sqrt{\alpha^2-\beta^2}\,(b-a)^{r+1/2}
           \bigg(\int_a^b\big|f^{(r)}(x)\big|^2\bigg)^{1/2}\mathrm dx\\
      &&\qquad\ge\,c_r\sqrt{\alpha^2-\beta^2}\,
           \bigg(\int_a^b\big|f^{(r)}(x)\big|^{1/(r+1)}\mathrm dx\bigg)^{r+1},
\end{eqnarray*}
where the first inequality follows from the Schwarz inequality and the second one from H\"older's 
inequality for integrals. As shown in the introduction, the gain can be significant, especially 
when the derivative $f^{(r)}$ drastically changes. 

\medskip
The approximation $\overline M^{\,*}_{N,r}$ possesses good asymptotic properties, but is not feasible 
since we do not have a direct access to the $C_i$s. In a feasible implementation one can approximate 
$C_i$ using divided differences, i.e.,
\begin{equation}\label{Citilda}
 \widetilde C_i=h\sqrt{\alpha^2-\beta^2}\,|d_i|\,r!\,\qquad\mbox{where}\qquad
d_i=f[x_{i,0},x_{i,1},\ldots,x_{i,r}]\end{equation}
and $x_{i,j}$ are arbitrary points from $I_i.$ Then 
$$ N_i^*=\frac{|d_i|^{1/(r+1)}}{\sum_{j=1}^{k_N}|d_j|^{1/(r+1)}}\,N.$$
This works well for functions $f$ for which the $r$th derivative does not nullify at any point in $[a,b].$ 
Indeed, then $f^{(r)}$ does not change its sign and, moreover, it is separated away from zero. 
This means that
$$\lim_{N\to\infty}\,\max_{1\le i\le k_N}\,{C_i}/{\widetilde C_i}=1,$$
which is enough for the asymptotic equality of Theorem \ref{thm:strata} to hold true. 

If $f^{(r)}$ is not separated away from zero then we may have problems with proper approximations 
of $C_i$ in the intervals $I_i$ where $|f^{(r)}|$ assumes extremely small values or even zeros. A possible and 
simple remedy is to choose `small' $\Delta>0$ and modify $\widetilde C_i$ as follows:
\begin{equation}\label{citil}\widetilde C_i=\left\{\begin{array}{rl}
  h\,\sqrt{\alpha^2-\beta^2}\,|d_i|\,r! &\;\,\mbox{for}\;|d_i|r!\ge\Delta,\\
  h\,\alpha\,\Delta\,r! &\;\,\mbox{for}\;|d_i|r!<\Delta.\end{array}\right.\end{equation} 
Then, letting $A_1=\big\{a\le x\le b:\,|f^{(r)}(x)|\ge\Delta\big\}$ and $A_2=[a,b]\setminus A_1,$ 
we have as $k\to+\infty$ that
\begin{eqnarray*}
 \bigg(\sum_{i=1}^kC_i^{1/(r+1)}\bigg)^{r+1} &\lessapprox& 
 \bigg(\sum_{i=1}^k\widetilde C_i^{1/(r+1)}\bigg)^{r+1} \\ &\approx&
    h^{-r}\,(\alpha^2-\beta^2)^{1/2}\bigg(\int_{A_1}\big|f^{(r)}(x)\big|^{\frac{1}{r+1}}\mathrm dx+
    |A_2|\Big(\sqrt{\tfrac{\alpha^2}{\alpha^2-\beta^2}}\,\Delta\Big)^{\frac{1}{r+1}}\bigg)^{r+1}.
\end{eqnarray*} 
Hence, the approximation of $C_i$ by \eqref{citil} results in an algorithm whose error is approximately 
upper bounded by 
$$ c_r\,\sqrt{\alpha^2-\beta^2}\,
    \bigg(\int_a^b\big|f_\Delta^{(r)}(x)\big|^{1/(r+1)}\mathrm dx\bigg)^{r+1}N^{-(r+1/2)},$$ where 
\begin{equation}\label{efdelta}
\big|f_\Delta^{(r)}(x)\big|=\max\bigg(\big|f^{(r)}(x)\big|,
\sqrt{\tfrac{\alpha^2}{\alpha^2-\beta^2}}\,\Delta\bigg).
\end{equation}
We obviously have $\lim_{\Delta\to 0^+}\int_a^b\big|f_\Delta^{(r)}(x)\big|^{1/(r+1)}\mathrm dx
=\int_a^b\big|f^{(r)}(x)\big|^{1/(r+1)}\mathrm dx.$

\medskip
A closer look at the deterministic part of $\overline M^{\,*}_{N,r}$ shows that the final 
partition of the interval $[a,b]$ tends to equalize the $L^1$ errors in all of the $m$ subintervals. 
As shown in \cite{PlaskotaSamoraj2022}, such a partition is optimal in the sense that it minimizes 
the asymptotic constant in the error $\|f-L_{m,r}f\|_{L^1(a,b)}$ among all possible piecewise 
Langrange interpolations $L_{m,r}$. A disadvantage of the optimal partition is that it is not nested. 
This makes it necessary to start the computations from scratch when $N$ is updated to a higher value. 
Also, a proper choice of the sequence $k_N=N^\kappa$ is problematic, especially when $N$ is still 
relatively small. On one hand, the larger $\kappa$ the better the approximation of $C_i$ by 
$\widetilde C_i,$ but also the more far away the partition from the optimal one. On the other hand, 
the smaller $\kappa$ the closer the partition to the optimal one, but also the worse the approximation 
of $C_i.$ This trade--off significantly affects the actual behavior of the algorithm, which can be seen 
in numerical experiments of Section \ref{sec:impl}. 

In the following section, we propose another adaptive approach leading to an easily implementable algorithm 
that produces nested partitions close to optimal and possesses asymptotic properties similar to that 
of $\overline M^{\,*}_{N,r}.$ As we shall see in Section \ref{sec:automatic}, nested partitions are vital for 
automatic Monte Carlo integration.  

\section{Second adaptive algorithm}\label{sec:second}

Consider a $\varrho$-weighted integration problem
$$S_\varrho f=\int_a^bf(x)\varrho(x)\,\mathrm dx,$$
where the function $\varrho:[a,b]\to\mathbb{R}$ is integrable and positive 
a.e. and $\int_a^b\varrho(x)\,\mathrm dx=1.$ The corresponding Monte Carlo algorithm is
$$M_{n,\varrho}f=\frac1n\sum_{i=1}^nf(t_i),\qquad t_i\stackrel{iid}\sim\mu_\varrho,$$
where $\mu_\varrho$ is the probability distribution on $[a,b]$ with density $\varrho.$ Then 
$$\mathbb E(S_\varrho f-M_{n,\varrho}f)^2=\frac1n\big(S_\varrho(f^2)-(S_\varrho f)^2\big).$$

Now, the non-weighted integral \eqref{theproblem} can be written as
$$S(f)=\int_a^b h(x)\varrho(x)\,\mathrm dx=S_\varrho(h),
\quad\mbox{where}\quad h(x)=\frac{f(x)}{\varrho(x)}.$$
Then
$$\mathbb E(Sf-M_{n,\varrho}h)^2=\mathbb E(S_\varrho h-M_{n,\varrho}h)^2
=\frac1n\big(S_\varrho(h^2)-(S_\varrho h)^2\big)=\frac1n\big(S_\varrho(f/\varrho)^2-(Sf)^2\big).$$ 
Let's go further on and apply a variance reduction, 
\begin{equation}\label{VRrho}
\overline M_{n,\varrho}f=S(Lf)+M_{n,\varrho}\bigg(\frac{f-Lf}{\varrho}\bigg),
\end{equation}
where $Lf$ is an approximation to $f.$ Then 
$$\mathbb E\big(Sf-\overline M_{n,\varrho}f\big)^2=
\frac1n\left(\int_a^b\frac{(f-Lf)^2(x)}{\varrho(x)}\,\mathrm dx-
\bigg(\int_a^b (f-Lf)(x)\,\mathrm dx\bigg)^2\right).$$
The question is how to choose $L$ and $\varrho$ to make the quantity 
$$\int_a^b\frac{(f-Lf)^2(x)}{\varrho(x)}\,\mathrm dx-
\bigg(\int_a^b (f-Lf)(x)\,\mathrm dx\bigg)^2$$
as small as possible.

Observe that if $$\varrho(x)=\frac{|(f-Lf)(x)|}{\|f-Lf\|_{L^1(a,b)}}$$ then 
$$\mathbb E\big(Sf-\overline M_{n,\varrho}f\big)^2=
\frac1n\left(\|f-Lf\|_{L^1(a,b)}^2-\bigg(\int_a^b(f-Lf)(x)\mathrm dx\bigg)^2\right)$$
and this error is even zero if $(f-Lf)(x)$ does not change its sign. This suggests the following algorithm.

Suppose that $Lf=L_{m,r}f$ is based on a partition of $[a,b]$ such that the $L^1$ errors 
in all $m$ subintervals $I_i$ have the same value, i.e.,
\begin{equation}\label{optpartition}
\|f-L_{m,r}f\|_{L^1(I_i)}=\frac1m\,\|f-L_{m,r}f\|_{L^1(a,b)},\quad 1\le i\le m.
\end{equation}
Then we apply the variance reduction \eqref{VRrho} with density
\begin{equation}\label{eq:rho}
\varrho(x)=\frac1{mh_i},\qquad x\in I_i,\quad1\le i\le m,
\end{equation}
where $h_i$ is the length of $I_i.$ That is, for the corresponding probability measure $\mu_\varrho$ 
we have $\mu_\varrho(I_i)=\tfrac1m$ and the conditional distribution $\mu_\varrho(\cdot|I_i)$ is uniform 
on $I_i.$ 

We now derive an error formula for such an approximation. 
Let $\gamma=\|P\|_{L^1(a,b)}=\int_0^1|P(z)|\,\mathrm dz.$ 
(Recall that $P$ is given by \eqref{mainpoly}.) We have
$$\|f-L_{m,r}f\|_{L^1(I_i)}=\frac{\gamma}{r!}\,h_i^{r+1}\big|f^{(r)}(\xi_i)\big|\quad\mbox{and}\quad
\|f-L_{m,r}f\|_{L^2(I_i)}=\frac{\alpha}{r!}\,h_i^{r+1/2}\big|f^{(r)}(\zeta_i)\big|$$
for some $\xi_i,\zeta_i\in I_i.$ Denoting
\begin{equation}\label{xidef} A=h_i^{r+1}\big|f^{(r)}(\xi_i)\big|\end{equation}
(which is independent of $i$) we have as $m\to+\infty$ that
\begin{eqnarray*}
\lefteqn{\bigg(\int_a^b\frac{(f-L_{m,r}f)^2(x)}{\varrho(x)}\,\mathrm dx\bigg)^{1/2}\;=\;
\bigg(m\sum_{i=1}^mh_i\int_{I_i}(f-L_{m,r}f)^2(x)\,\mathrm dx\bigg)^{1/2}} \\
&&\;=\;\frac{\alpha}{r!}\,\bigg(m\sum_{i=1}^m h_i^{2r+2}\big|f^{(r)}(\zeta_i)\big|^{2}\bigg)^{1/2}
\,\approx\,\frac{\alpha}{r!}\,\bigg(m\sum_{i=1}^m h_i^{2r+2}\big|f^{(r)}(\xi_i)\big|^2\bigg)^{1/2} \\
&&\;=\;\frac{\alpha}{r!}\,\bigg(m\sum_{i=1}^mA^2\bigg)^{1/2}\,=\,\frac{\alpha}{r!}\,mA\,=\,
\frac{\alpha}{r!}\,\big(mA^{1/(r+1)}\big)^{r+1}m^{-r} \\
&&\;=\;\frac{\alpha}{r!}\,\bigg(\sum_{i=1}^m h_i\big|f^{(r)}(\xi_i)\big|^{1/(r+1)}\bigg)^{r+1}m^{-r}
\;\approx\;\frac{\alpha}{r!}\,\bigg(\int_a^b\big|f^{(r)}(x)\big|^{1/(r+1)}\mathrm dx\bigg)^{r+1}m^{-r}.
\end{eqnarray*}

To get an asymptotic formula for $\int_a^b(f-L_{m,r}f)(x)\,\mathrm dx$ we use the analysis done 
in Section \ref{sec:equispaced}. If $\beta=0$ then the integral decreases faster than $m^{-r}.$ 
Let $\beta\ne 0.$ Then 
$$\int_a^b(f-L_{m,r}f)(x)\,\mathrm dx\;\approx\;
\frac{\beta}{r!}\,\sum_{i=1}^m h_i^{r+1}f^{(r)}(\xi_i)=\frac{\beta}{r!}(m_+-m_-)A,$$
where $\xi_i$s are as in \eqref{xidef}, and $m_+$ and $m_-$ are the numbers of indexes $i$ for which 
$f^{(r)}(\xi_i)\ge 0$ and $f^{(r)}(\xi_i)<0,$ respectively. Let 
$$D_+=\{x\in[a,b]:\,f^{(r)}(x)\ge 0\},\qquad D_-=\{x\in[a,b]:\,f^{(r)}(x)<0\}.$$
Since $A\approx\|f^{(r)}\|_{L^{1/(r+1)}(a,b)}m^{-(r+1)}$ and
$m_+A^{1/(r+1)}\approx\int_{D_+}|f^{(r)}(x)|^{1/(r+1)}\mathrm dx,$ we have 
$$m_+\approx\frac{\int_{D_+}|f^{(r)}(x)|^{1/(r+1)}\mathrm dx}
   {\int_{D_+\cup D_-}|f^{(r)}(x)|^{1/(r+1)}\mathrm dx},\qquad 
     m_-\approx\frac{\int_{D_-}|f^{(r)}(x)|^{1/(r+1)}\mathrm dx}
  {\int_{D_+\cup D_-}|f^{(r)}(x)|^{1/(r+1)}\mathrm dx}.$$
Thus
$$\int_a^b(f-L_{m,r}f)(x)\,\mathrm dx\approx\frac{\beta}{r!}
\bigg(\frac{\int_a^b|f^{(r)}(x)|^{1/(r+1)}\mathrm{sgn} f^{(r)}(x)\,\mathrm dx}
{\int_a^b|f^{(r)}(x)|^{1/(r+1)}\mathrm dx}\bigg)\|f^{(r)}\|_{L^{1/(r+1)}(a,b)}\,m^{-r}$$
provided $\int_a^b|f^{(r)}(x)|^{1/(r+1)}\mathrm{sgn} f^{(r)}(x)\,\mathrm dx\ne 0,$ 
and otherwise the convergence is faster than $m^{-r}.$

\smallskip
Our analysis above shows that if $m$ and $n$ are chosen as in 
\eqref{mn1} and \eqref{mn2} then the error of the described algorithm asymptotically 
equals (as $N\to+\infty$)
\begin{equation}\label{annoy}
 c_r\,\sqrt{\alpha^2-\beta^2
\bigg(\frac{\int_a^b|f^{(r)}(x)|^{1/(r+1)}\mathrm{sgn}f^{(r)}(x)\,\mathrm dx}
{\int_a^b|f^{(r)}(x)|^{1/(r+1)}\mathrm dx}\bigg)^2}
\bigg(\int_a^b\big|f^{(r)}(x)\big|^{1/(r+1)}\mathrm dx\bigg)^{r+1}N^{-(r+1/2)}. 
\end{equation}

\smallskip
The factor at $\beta^2$ in \eqref{annoy} can be easily replaced by $1$ with the help 
of stratified sampling. Indeed, instead of randomly sampling $n$ times with density 
\eqref{eq:rho} on the whole interval $[a,b],$ one can apply the same sampling strategy 
independently on $k$ groups $G_j$ of subintervals. Each group consists of $s=m/k$ 
subintervals, 
$$G_j=\bigcup_{\ell=1}^s I_{(j-1)s+\ell},\qquad 1\le j\le k,$$
and the number of samples for each $G_j$ equals $n/k.$ As in the algorithm 
$\overline M^{\,*}_{N,r},$ we combine $k=k_N$ and $N$ in such a way that 
$k_N\to+\infty$ and $k_N/N\to 0$ as $N\to\infty.$ 
Then the total number of points used in each $G_j$ is $N_j=N/k.$ Denoting
$$C_j=\bigg(\int_{G_j}\big|f^{(r)}(x)\big|^{1/(r+1)}\mathrm dx\bigg)^{r+1}
=\bigg(\frac 1k\int_a^b\big|f^{(r)}(x)\big|^{1/(r+1)}\mathrm dx\bigg)^{r+1}$$
and using the fact that the factor at $\beta^2$ equals $1$ if $f^{(r)}$ does not 
change its sign, the error of such an approximation asymptotically equals
$$c_r\,\sqrt{\alpha^2-\beta^2}\,\bigg(\sum_{j=1}^k\frac{C_j^2}{N_j^{2r+1}}\bigg)^{1/2}
=c_r\,\sqrt{\alpha^2-\beta^2}\,
\bigg(\int_a^b\big|f^{(r)}(x)\big|^{1/(r+1)}\mathrm dx\bigg)^{r+1}N^{-(r+1/2)},$$
as claimed. 
(Note that $N_j=N/k$ minimize the sum $\sum_{j=1}^kC_j^2N_j^{-(2r+1)}$ 
with respect to $\sum_{j=1}^kN_j=N;$ compare with the analysis 
in Section \ref{sec:strata}.)

Thus we obtained exactly the same error formula as in 
Theorem \ref{thm:strata} for $\overline M_{N,r}^{\,*}.$

\medskip
It remains to show a feasible construction of a nested partition that is close to the one 
satisfying \eqref{optpartition}. To that end, we utilize the iterative method presented in 
\cite{PlaskotaSamoraj2022}, where the $L^p$ error of piecewise Lagrange interpolation is examined. 

\smallskip
We first consider the case when
\begin{equation}\label{derpos}
f^{(r)}>0\quad\mbox{or}\quad f^{(r)}<0.
\end{equation}
In the following construction, we use a priority queue $\mathcal S$ whose elements are subintervals. 
For each subinterval $I_i$ of length $h_i,$ its priority is given as
\begin{equation*}
p_f(I_i)\,=\,h_i^{r+1}|d_i|, 
\end{equation*}
where $d_i$ is the divided difference \eqref{Citilda}. 
In the following pseudocode, $\mathrm{insert}(\mathcal S,I)$ and 
$I:=\mathrm{extract\_max}(\mathcal S)$ implement correspondingly 
the actions of inserting an interval to $\mathcal S,$ and extracting 
from $\mathcal S$ an interval with the highest priority.

\medskip
$\mathbf{algorithm}\;\mathrm{PARTITION}$

$\mathcal S=\emptyset;\;\mathrm{insert}(\mathcal S,[a,b]);$ 

$\mathbf{for}\; k=2:m$ 

$\quad [l,r]=\mathrm{extract\_max}(\mathcal S);$ 

$\quad c=(l+r)/2;$ 

$\quad \mathrm{insert}(\mathcal S,[l,c]); \mathrm{insert}(\mathcal S,[c,r]);$ 

$\mathbf{endfor}$

\medskip\noindent
After execution, the elements of $\mathcal S$ form a partition into $m$ subintervals $I_i.$ 
Note that if the priority queue is implemented through a \emph{heap} then the running time
of $\mathrm{PARTITION}$ is proportional to $m\log m.$

\smallskip
Denote by $\overline M^{\,**}_{N,r}$ the corresponding algorithm that 
uses the above nested partition and density \eqref{eq:rho}, and $N$ is related to the number 
$m$ of subintervals and the number $n$ of random samples as in \eqref{mn1} and \eqref{mn2}. 
We want to see how much worse is this algorithm than that using the (not nested) partition 
\eqref{optpartition}.

Let $A=(A_1,A_2,\ldots,A_m)$ with 
$$A_i=p_f(I_i)\,r!=h_i^{r+1}\big|f^{(r)}(\omega_i)\big|,\qquad\omega_i\in I_i,$$
and $\|A\|_p=\big(\sum_{i=1}^mA_i^p\big)^{1/p}.$ For the corresponding piecewise Lagrange 
approximation $L_{m,r}f$ and density $\varrho$ given by \eqref{eq:rho} we have 
\begin{eqnarray*}
\lefteqn{\int_a^b\frac{(f-L_{m,r}f)^2(x)}{\varrho(x)}\,\mathrm dx
\,-\,\bigg(\int_a^b(f-L_{m,r}f)(x)\,\mathrm dx\bigg)^2}\\
&&\approx\,\frac{1}{(r!)^2}\bigg(\alpha^2 m\sum_{i=1}^m A_i^2
-\beta^2\bigg(\sum_{i=1}^mA_i\bigg)^2\,\bigg)\,=\,\frac{1}{(r!)^2}
\left(\,\alpha^2m\|A\|_2^2-\beta^2\|A\|_1^2\,\right).
\end{eqnarray*}
We also have 
$\bigg(\int_a^b\big|f^{(r)}(x)\big|^{1/(r+1)}\mathrm dx\bigg)^{r+1}\approx
\big(\sum_{i=1}^mA_i^{1/(r+1)}\big)^{r+1}=\|A\|_{\frac{1}{r+1}}.$ Hence
\begin{eqnarray*}\nonumber
\sqrt{\mathbb E(Sf-\overline M^{**}_{N,r}f)^2}&\approx&
\left(\,\alpha^2m\|A\|_2^2-\beta^2\|A\|_1^2\,\right)^{1/2}n^{-1/2}m^{-r}\\
&\approx&K_{m,r}(A)\,c_r\,\sqrt{\alpha^2-\beta^2}\,
\bigg(\int_a^b\big|f^{(r)}(x)\big|^{1/(r+1)}\bigg)^{r+1}N^{-(r+1/2)},
\end{eqnarray*}
where $$K_{m,r}(A)=\frac{\sqrt{\kappa_\alpha^2\,m\,\|A\|_2^2
-\kappa_\beta^2\,\|A\|_1^2}}{\|A\|_{\frac{1}{r+1}}}\,m^r,
\qquad\kappa_\alpha=\frac{\alpha}{\sqrt{\alpha^2-\beta^2}},
\quad\kappa_\beta=\frac{\beta}{\sqrt{\alpha^2-\beta^2}}.$$ 

Observe that halving an interval results in two subintervals whose priorities are asymptotically 
(as $m\to+\infty$)  $2^{r+1}$ times smaller than the priority of the original interval. 
This means that $K_{m,r}(A)$ is asymptotically not larger than
\begin{equation}\label{KA}
K^*(r)\;=\;\limsup_{m\to\infty}\;\max\,
\left\{K_{m,r}(A):\;A=(A_1,\ldots,A_m),\,\max_{1\le i,j\le m}\frac{A_i}{A_j}\le 2^{r+1}\,\right\}.
\end{equation}
Thus we obtained the following result.

\begin{theorem}\label{thm:importstrata}
If the function $f$ satisfies \eqref{derpos} then we have as $N\to+\infty$ that
$$\sqrt{\mathbb E(Sf-\overline M^{\,**}_{N,r}f)^2}\,\lessapprox\,K^*(r)\,c_r\,\sqrt{\alpha^2-\beta^2}\,
\bigg(\int_a^b\big|f^{(r)}(x)\big|^{1/(r+1)}\mathrm dx\bigg)^{r+1}N^{-(r+1/2)},$$
where $K^*(r)$ is given by \eqref{KA}.
\end{theorem}

We numerically calculated $K^*(r)$ in some special cases. 
For instance, if the points $z_i$ in \eqref{zi-points} are equispaced, $z_i=(i-1)/(r-1),$ $1\le i\le r,$ 
then for $r=2,3,4,5,6$ we correspondingly have 
$$K^*(r)\,=\,4.250,\;3.587,\;7.077,\;11.463,\;23.130,$$ 
while for any $z_i$s satisfying $\beta=\int_0^1(z-z_1)\cdots(z-z_r)\,\mathrm dz=0$ we have 
$$K^*(r)\,=\,2.138,\;3.587,\;6.323,\;11.463,\;21.140.$$

\medskip
If $f$ does not satisfy \eqref{derpos} then the algorithm $\overline M^{**}_{N,r}f$ may fail. 
Indeed, it may happen that $p_f(I_i)=0$ while $f^{(r)}\not=0$ in $I_i.$ Then this subinterval may 
never be further subdivided. In this case, we can repeat the same construction, 
but with the modified priority
$$p_f(I_i)\,=\,h_i^{r+1}\max\big(\,|d_i|,\,\Delta/r!\big),$$
where $\Delta>0.$ Then the error is asymptotically 
upper bounded by
$$K^*(r)\,c_r\,\sqrt{\alpha^2-\beta^2}\,
\bigg(\int_a^b\big|f_\Delta^{(r)}(x)\big|^{1/(r+1)}\mathrm dx\bigg)^{r+1}N^{-(r+1/2)},$$
where $\big|f_\Delta^{(r)}(x)\big|$ is given by \eqref{efdelta}.

\section{Numerical experiments}\label{sec:impl}

In this section, we present results of two numerical experiments that illustrate the performance 
of the nonadaptive Monte Carlo algorithm $\overline M_{N,r}$ and the adaptive algorithms 
$\overline M^{\,*}_{N,r}$ and $\overline M^{**}_{N,r}.$ Our test integral is 
$$\int_0^1\frac{1}{x+10^{-4}}\,\mathrm dx.$$
Since for $r\in\mathbb{N}$ we have $(-1)^r f^{(r)}>0,$ the parameter $\Delta$ is set to zero.

The three algorithms are verified for $r=2$ and $r=4$. In both cases, the interpolation nodes
are equispaced, i.e., in \eqref{zi-points} we take $$z_i =\frac{i-1}{r-1},\qquad 1\leq i \leq r.$$ 
In addition, for the first adaptive algorithm $\overline M_{N,r}^{\,*}$ we take $k_N = N^\kappa$ 
with $\kappa = 0.8.$ This exponent was chosen to ensure a trade--off as per our discussion 
in Section \ref{sec:strata}, and some empirical results. Also, for a fixed $N$ we plot a single output 
instead of the expected value estimator. Therefore the error fluctuations are visible. For completeness, 
we also show the asymptotes corresponding to the theoretical errors from Theorems \ref{thm:equ} 
and \ref{thm:strata}, and the upper bound from Theorem \ref{thm:importstrata}. The scale is 
logarithmic, $-\log_{10}(\mathrm{error})$ versus $\log_{10}N.$ 

\begin{figure}[!htp]
    \centering
    \includegraphics[width=0.7\textwidth]{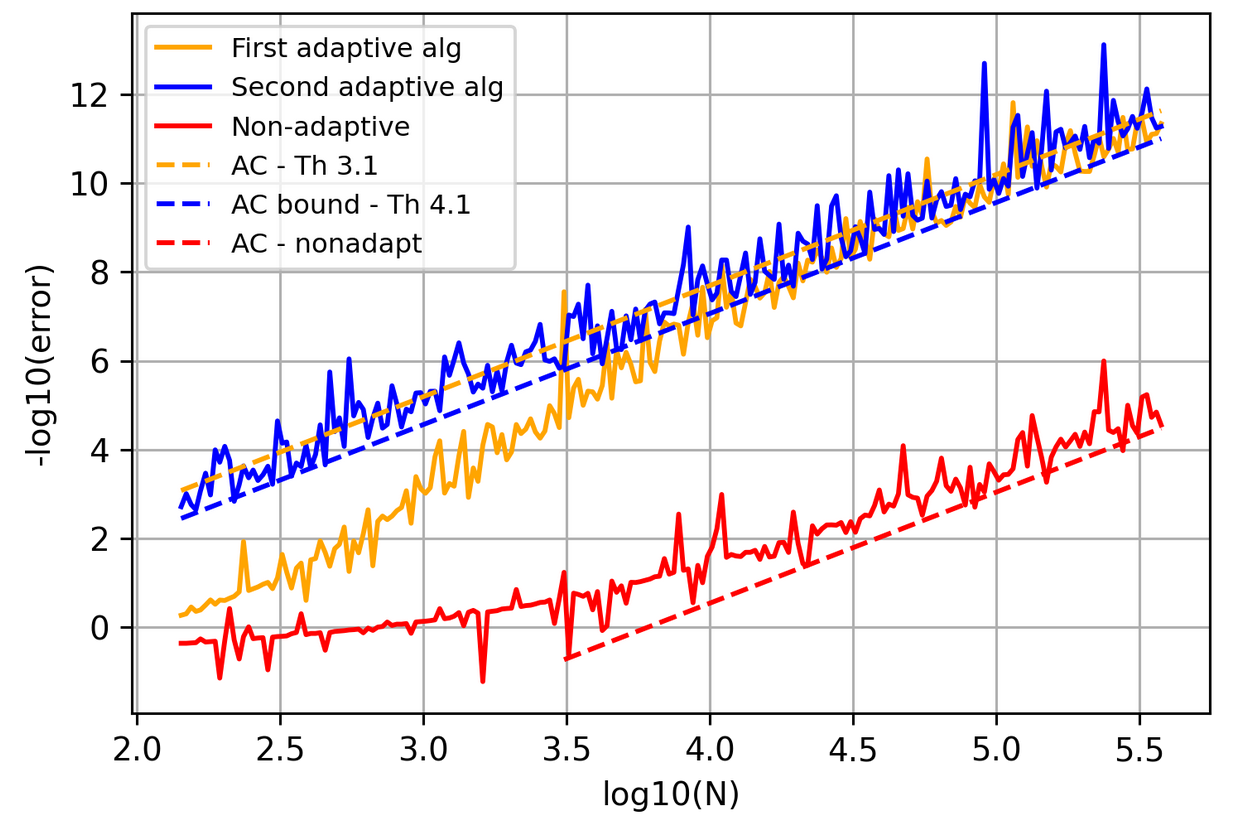}
    \caption{Comparison of nonadaptive and adaptive Monte Carlo algorithms together with 
    related asymptotic constants (AC) for $r=2.$}
    \label{impl:MC_test_r2}
\end{figure}

The results for $r=2$ are presented in Figure~\ref{impl:MC_test_r2}. 
\FloatBarrier

As it can be observed, both adaptive algorithms significantly outperform the nonadaptive MC; however, 
the right asymptotic behaviour of the first adaptive algorithm is visible only for large $N.$ 

Similar conclusions can be inferred from validation performed for $r=4,$ with all other parameters unchanged. 
We add that the results for $N$ larger than $10^{4.8}$ are not illustrative any more, since the process is 
disturbed by a serious reduction of significant digits when calculating divided differences in the partition part.

\begin{figure}[!htp]
    \centering
    \includegraphics[width=0.7\textwidth]{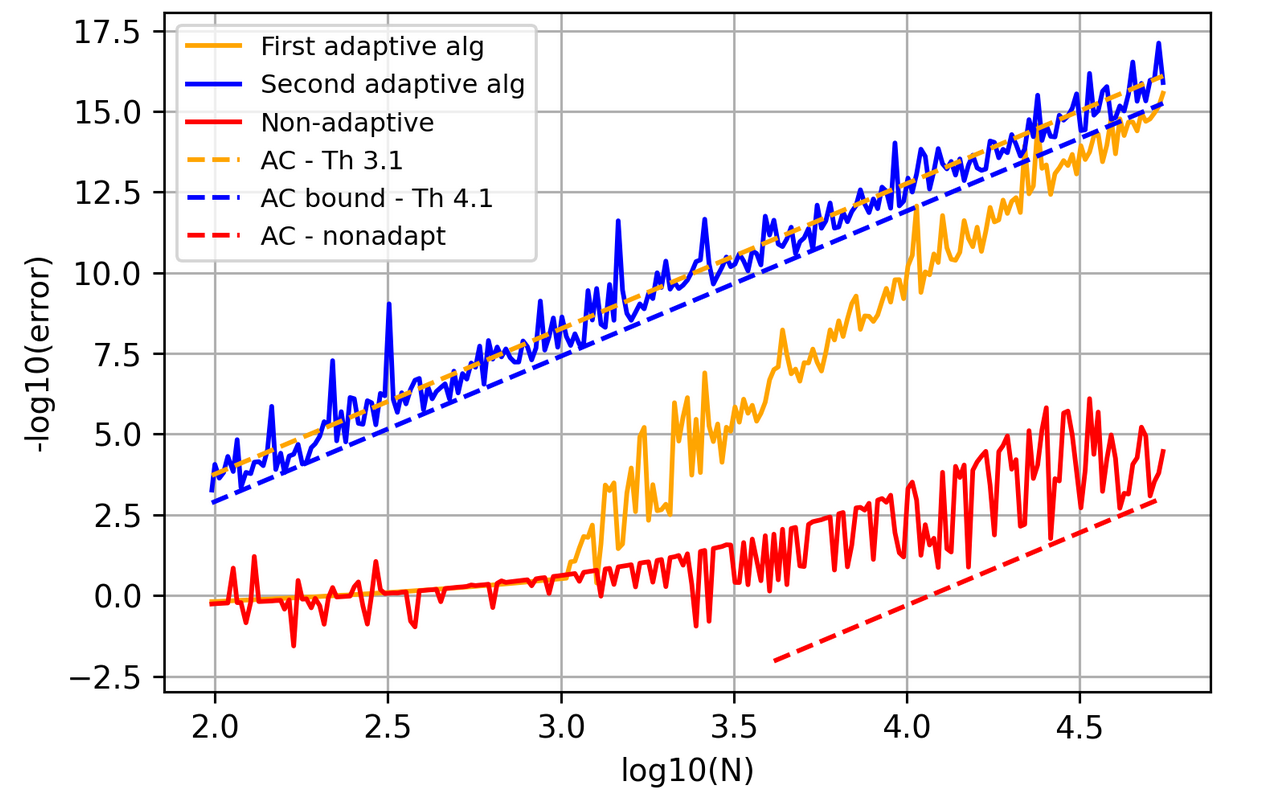}
    \caption{Comparison of nonadaptive and adaptive Monte Carlo algorithms together with related 
    asymptotic constants (AC) for $r=4.$}
    \label{impl:MC_test_r4}
\end{figure}
\FloatBarrier
Notably, both adaptive algorithms attain their asymptotic errors, but this is not the case for nonadaptive MC 
for which the output is not stable. Initially, the first adaptive algorithm does not leverage additional sampling since 
for all intervals $I_i$ we have $N_i/(2r+1)< 1.$ The Monte Carlo adjustments are visible only for $N \geq 10^3$ 
and the error tends to the theoretical asymptote.

In conclusion, the numerical experiments confirm our theoretical findings and, in particular, superiority of 
the second adaptive algorithm $\overline M^{\,**}_{N,r}.$

\section{Automatic Monte Carlo integration}\label{sec:automatic}

We now use the results of Section \ref{sec:second} for automatic Monte Carlo integration. 
The goal is to construct an algorithm that for given $\varepsilon>0$ and $0<\delta<1$ returns 
an $\varepsilon$-approximation to the integral $Sf$ with probability at least $1-\delta,$ 
asymptotically as $\varepsilon\to 0^+.$ To that end, we shall use the approximation 
$\overline M_{N,r}^{\,**}f$ with $N$ determined automatically depending on $\varepsilon$ and $\delta.$ 

\smallskip
Let $X_i$ for $1\le i\le n$ be independent copies of the random variable
$$X=S(f-L_{m,r}f)-\frac{(f-L_{m,r}f)(t)}{\varrho(t)},\qquad t\sim\mu_{\varrho},$$
where $L_{m,r},$ $n$ and $\varrho$ are as in $\overline M_{N,r}^{\,**}f.$ 
Then $\mathbb E(X)=0$ and
$$Sf-\overline M_{N,r}^{\,**}f=\frac{X_1+X_2+\cdots+X_n}{n}.$$
By Hoeffding's inequality \cite{Hoeffding1963} we have
$$\mathrm{Prob}\left(\big|Sf-\overline M_{N,r}^{\,**}f\big|>\varepsilon\right)\,\le\,
2\,\exp\left(\frac{-\varepsilon^2n}{2\,B_m^2}\right),$$
where $B_m=\max_{a\le t\le b}|X(t)|.$ Hence we fail with probability at most $\delta$ if 
\begin{equation}\label{Hoeff}
\frac{\varepsilon^2n}{2B_m^2}\,\ge\,\ln\frac{2}{\delta}.
\end{equation}

Now we estimate $B_m.$ Let $\lambda=\|P\|_{L^\infty(0,1)}=\max_{0\le t\le 1}|P(t)|,$ and 
$$\mathcal L_r(f)=\bigg(\int_a^b\big|f^{(r)}_\Delta(x)\big|^{1/(r+1)}\mathrm dx\bigg)^{r+1},$$ 
where $\Delta=0$ if $f^{(r)}>0$ or $f^{(r)}<0,$ and $\Delta>0$ otherwise. Let 
$A=(A_1,A_2,\ldots,A_m)$ with
$$A_i=h_i^{r+1}\max_{x\in I_i}|f^{(r)}(x)|,\quad 1\le i\le m,$$
where, as before, $\{I_i\}_{i=1}^m,$ is the partition used by $\overline M_{N,r}^{\,**}f$ 
and $h_i$ is the length of $I_i.$ Since 
$\|A\|_\frac{1}{r+1}=\big(\sum_{i=1}^mA_i^{1/(r+1)}\big)^{r+1}\lessapprox\mathcal L_r(f),\;$ 
for $x\in I_i$ we have  
$$\frac{\big|f(x)-L_{m,r}f(x)\big|}{\varrho(x)}\,\le\,\frac{\lambda}{r!}\,m\,A_i\,\lessapprox\,
\frac{\lambda}{r!}\,\bigg(\frac{m\,\|A\|_\infty}{\|A\|_\frac{1}{r+1}}\bigg)\mathcal L_r(f)
\,\lessapprox\,2^{r+1}\frac{\lambda}{r!}\,\mathcal L_r(f)\,m^{-r}.$$ 
We have the same upper bound for $S(f-L_{m,r}f)$ since by mean-value theorem 
$$S(f-L_{m,r}f)=\int_a^b\frac{(f-L_{m,r}f)(x)}{\varrho(x)}\,\varrho(x)\,\mathrm dx
=\frac{(f-L_{m,r}f)(\xi)}{\varrho(\xi)},\qquad\xi\in[a,b].$$
Hence
$$B_m\,\lessapprox\,2^{r+2}\frac{\lambda}{r!}\,\mathcal L_r(f)\,m^{-r}.$$

Using the above inequality and the fact that $\sqrt{n}\,m^r\approx N^{r+1/2}/(c_rr!)$ 
with $c_r$ given by \eqref{ciar}, we get 
$$\frac{\varepsilon^2n}{2B_m^2}\gtrapprox\left(\frac{\varepsilon\,N^{r+1/2}}
{\hat c_r\,\mathcal L_r(f)}\right)^{\!\!2},\quad\mbox{where}\quad
\hat c_r=2^{r+5/2}\lambda c_r.$$
The last inequality and \eqref{Hoeff} imply that we fail to have error $\varepsilon$ with probability at most $\delta$ for
\begin{equation}\label{eNform}
N\,\gtrapprox\,\left(\hat c_r\,\mathcal L_r(f)\,\frac{\sqrt{\ln(2/\delta)}}{\varepsilon}\right)^{\frac1{r+1/2}},
\qquad\mbox{as}\quad\varepsilon\to 0^+.
\end{equation}

\smallskip
Now the question is how to obtain the random approximation $\overline M_{N,r}^{\,**}f$ for $N$ satisfying \eqref{eNform}. 

\medskip
One possibility is as follows. We first execute the iteration $\mathbf{for}$ in the algorithm $\mathrm{PARTITION}$ 
of Section \ref{sec:second} for $k=2:m,$ where $m$ satisfies 
$\lim_{\varepsilon\to 0^+}m\,\varepsilon^{\frac{1}{r+1/2}}=0,$ e.g.,
$$m=\left\lfloor\bigg(\frac{\sqrt{\ln(2/\delta)}}{\varepsilon}\bigg)^\frac{1}{r+1}\right\rfloor.$$
Let $\{I_i\}_{i=1}^m$ be the obtained partition. Then we replace $\mathcal L_r(f)$ in \eqref{eNform} 
by its asymptotic equivalent 
\begin{equation}\label{ElT}
\widetilde{\mathcal L}_r(f)=\bigg(\sum_{i=1}^{m}p_f(I_i)^\frac{1}{r+1}\bigg)^{r+1},
\end{equation} set 
\begin{equation}\label{eNeps}
N_\varepsilon=\left\lfloor\left(\hat c_r\,\widetilde{\mathcal L}_r(f)\,\frac{\sqrt{\ln(2/\delta)}}
{\varepsilon}\right)^{\frac1{r+1/2}}\right\rfloor,
\end{equation}
and continue the iteration for $k=m+1:m_\varepsilon,$
where $m_\varepsilon$ is the number of subintervals corresponding to $N_\varepsilon.$  
Finally, we complete the algorithm by $n_\varepsilon$ random samples.

Denote the final randomized approximation by $\mathcal A_{\varepsilon,\delta}f.$ 
Then we have $\mathcal A_{\varepsilon,\delta}f=\overline M_{N_\varepsilon,r}^{\,**}f$ and
$$\mathrm{Prob}\big(\,\big|Sf-\mathcal A_{\varepsilon,\delta}f|>\varepsilon\big)\,\lessapprox\,\delta,
\qquad\mbox{as}\quad\varepsilon\to 0^+.$$

\medskip
A disadvantage of the above algorithm is that it uses a priority queue and therefore its total running time 
is proportional to $N\log N.$ It turns out that by using recursion the running time can be reduced to $N.$ 
 
A crucial component of the algorithm with the running time 
proportional to $N$ is the following recursive procedure, in which $\mathcal S$ is a set of intervals.

\medskip
$\mathbf{procedure}\;\mathrm{AUTO}\,(f,a,b,e)$

$\mathbf{if}\;p_f([a,b])\le e$

$\quad \mathrm{insert}(\mathcal S,[a,b]);$

$\mathbf{else}$

$\quad c:=(a+b)/2;$

$\quad\mathrm{AUTO}(f,a,c,e);$

$\quad\mathrm{AUTO}(f,c,b,e);$

$\mathbf{endif}$

\medskip
Similarly to $\mathcal A_{\varepsilon,\delta},$ the algorithm consists of two steps. First $\mathrm{AUTO}$ 
is run for $e=\varepsilon'$ satisfying $\varepsilon'\to 0^+$ and $\varepsilon/\varepsilon'\to 0^+$ as $\varepsilon\to 0^+,$ e.g., 
$$\varepsilon'=\varepsilon^\kappa,
\quad\mbox{where}\quad 0<\kappa<1.$$ Then $\mathcal L_r(f)$ in \eqref{eNform} is replaced by 
$\widetilde{\mathcal L}_r(f)$ given by \eqref{ElT}, and $N_\varepsilon$ found from \eqref{eNeps}. 
The recursion is resumed with the target value $e=\varepsilon'',$ where
$$\varepsilon''=\widetilde{\mathcal L}_r(f)\,m_\varepsilon^{-(r+1)}.$$
The algorithm is complemented by the corresponding $n_\varepsilon$ random samples. 

Observe that the number $m''$ of subintervals in the final partition is asymptotically at least $m_\varepsilon.$ 
Indeed, for any function $g\in C^{r}([a,b])$ with $g^{(r)}(x)=\big|f_\Delta^{(r)}(x)\big|$ we have 
$\mathcal L_r(g)=\mathcal L_r(f)$ and 
\begin{eqnarray*}
\frac{\gamma}{r!}\,\widetilde{\mathcal L}_r(f)\big(m''\big)^{-r} &\lessapprox &
\|g-L_{m'',r}g\|_{L^1(a,b)}\,\approx\,\frac{\gamma}{r!}\sum_{i=1}^{m''}(h_i'')^{r+1}g^{(r)}(\xi_i) \\
&\approx &\frac{\gamma}{r!}\sum_{i=1}^{m''}p_f(I_i'')\lessapprox\frac{\gamma}{r!}\,m''\varepsilon'',
\end{eqnarray*}
where the first inequality above follows from Proposition 2 of \cite{PlaskotaSamoraj2022}. This implies
$$m''\,\gtrapprox\,\bigg(\frac{\widetilde{\mathcal L}_r(f)}{\varepsilon''}\bigg)^\frac{1}{r+1}\approx\,m_\varepsilon,$$
as claimed. 

Denote the resulting approximation by $\mathcal A_{\varepsilon,\delta}^*f.$ Observe that its running time 
is proportional to $N_\varepsilon$ since recursion can be implemented in linear time.

\begin{theorem}\label{thm:auto}
We have 
$$\mathrm{Prob}\big(\,\big|Sf-\mathcal A_{\varepsilon,\delta}^*f|>\varepsilon\big)\,\lessapprox\,\delta,
\qquad\mbox{as}\quad\varepsilon\to 0^+.$$
\end{theorem}

\bigskip
\phantomsection

Now we present outcomes of the second automatic procedure $\mathcal A_{\varepsilon,\delta}^*$ for 
the test integral
\begin{equation}\label{autoint}
\int_0^1\cos\bigg(\dfrac{100\,x}{x+10^{-4}}\bigg)\,\mathrm dx.
\end{equation}
Although the derivatives fluctuate and nullify many times in this case, we take $\Delta=~0.$
We confront the outcomes for $r=2$ and $r=4.$ In each case, we compute the number of breaches 
(i.e. when the absolute error is greater than $\varepsilon = 10^{-3}$) based on $K=10\,000$ independent 
executions of the code. We also take $\varepsilon' = \varepsilon^{1/2}.$ In our testing, we expect 
the empirical probability of the breach to be less than $\delta=0.05.$ For completeness, we also present 
the maximum error from all executions together with obtained $N_\varepsilon.$

\begin{table}
\caption{Performance of the second automatic algorithm for the integral \eqref{autoint}.}
\label{tab:automatic}
\begin{tabular}{llllll}
\hline\noalign{\smallskip}
 & $\varepsilon$ &  $\delta$ & $K$ & $N_\varepsilon$ & $e_{max}$ \\
\hline\noalign{\smallskip}
 $r = 2$ & $10^{-3}$ & 0.05 & 10\,000 & 3\,092 & $4.5 \cdot 10^{-5}$ \\
 $r = 4$ & $10^{-3}$ & 0.05 & 10\,000 & 811 & $5.5 \cdot 10^{-6}$ \\
 \hline
\end{tabular}
\end{table}
\FloatBarrier
Note that in both cases we did not identify any exceptions. The magnitude of the maximum errors indicate 
a serious overestimation of $N_\varepsilon,$ but the results are  satisfactory given the upper bound estimate 
of Theorem \ref{thm:auto}.

\section*{Acknowledgements}
The work of L. Plaskota and P. Przyby{\l}owicz 
was partially supported by the National Science Centre, Poland, under project 2017/25/B/ST1/00945.

\section*{Statements and Declarations}
All authors declare that they have no conflicts of interest.



\section*{Appendix}\label{sec:appendix}

Below we present a crucial part of the code in the Python programming language, where all the algorithms 
were implemented. In addition, we provide relevant comments linked to particular fragments of the code.

\begin{lstlisting}[language=Python, caption=Second adaptive algorithm $\overline M^{**}_{N,r}$- crucial 
part of the code. , basicstyle=\ttfamily\tiny]
def second_adaptive_MC(a, b, N, main_nodes, f, r, node_type = 'uniform'): #(1)
    partial_quad = Decimal('0.0')  #(2)
    
    if node_type == 'uniform':  #(3)
        n = int(np.floor((N-1)/(1 + 2 *r)))
        m = int(np.floor((2 * r * (N-1))/(2 * r ** 2 - r - 1)))
    
    mc_init = MC_samples_nonunif(m, n)  #(4)
    
    #(5)
    MnR = Decimal('0.0')    
    l = 0
    
    for i in range(len(main_nodes) - 1 ):
        h_i = main_nodes[i+1] - main_nodes[i]
        
        #(6)
        if node_type == 'uniform':
            interpol_ix = main_nodes[i] * np.ones(r) + np.multiply(optimalt_equidistant(r), h_i)
        
        interpol_iy = []
        for s in range(1,r+1):
            interpol_iy.append(f(interpol_ix[s-1]))
        
        #(7)
        if r == 2:
            SL_mr = Decimal('0.5') * Decimal(h_i) \ 
                    * (Decimal(f(main_nodes[i])) + Decimal(f(main_nodes[i+1])))
        elif r == 4:
            SL_mr = Decimal('0.125') * Decimal(h_i) * Decimal((f(main_nodes[i]) \
                  + 3 * f(main_nodes[i] + 1/3 * h_i) + 3 * f(main_nodes[i] + 2/3 * h_i) \ 
                  + f(main_nodes[i+1])))
        
        #(8)
        while l < n and mc_init[l] < (i+1):
            mc_point = math.modf(mc_init[l])[0] * h_i + main_nodes[i] #(9)
            MnR = MnR + (Decimal(f(mc_point)) - Decimal(lagrange(interpol_ix, interpol_iy, mc_point))) \
            * Decimal(h_i)  #(10)
            l = l + 1
        
        partial_quad = partial_quad + SL_mr #(11)
    
    #(12)
    MnR = Decimal(MnR) * Decimal(m) / Decimal(n)
    partial_quad = partial_quad + MnR
    return partial_quad
\end{lstlisting}

\begin{enumerate}
    \item [(1)] The almost optimal partition \texttt{main\_nodes} is derived out of this function in order 
    to save computation time when the trajectories are computed subsequently. Moreover, \texttt{node\_type} 
    argument lets the user insert his own partitions, e.g. those based on Chebyshev polynomials of the second kind.
    \item [(2)] In order to minimize errors resulting from (possibly) adding relatively small adjustments to 
    the estimated quadrature value, we use \texttt{Decimal} library. It enables us to increase the precision of
    intermediate computations, which is now set to 28 digits in decimal system.
    \item [(3)] In our case, the interpolating polynomial is based on equidistant mesh including endpoints of 
    a subinterval $I_i$. By \texttt{np} we understand the references to NumPy library.
    \item [(4)] Initializing the variables which control Monte Carlo adjustments for our quadrature. 
    In particular, \texttt{l} stores the number of currently used random points, while we loop through the subintervals. 
    \item [(5)] The program calculates all interpolation nodes in the interval $I_i.$ For that reason, the function 
    \texttt{optimalt\_uniform} is executed to provide distinct $z_1, \ldots, z_r \in [0,1].$
    \item [(6)] Depending on the value of $r,$ different formulas for (nonadaptive, deterministic) quadrature 
    $SL_{m,r}$ are leveraged.
    \item [(7)] Below, we calculate the Monte Carlo adjustment on interval $I_i.$
    \item [(8)] This code yields random points used for . \texttt{MC\_init} function reports them in a from of 
    a number from 0 to $m.$ The integer part points the index $i$ of subinterval, while the fractional part - its 
    position within $I_i.$ Both parameters are sourced by using \texttt{math.modf} function.  
    \item [(9)] For stability reasons, the coefficients of interpolating polynomial in canonical base are not stored. 
    Therefore, for every point, \texttt{lagrange} function is invoked separately.
    \item [(10)] We decided to add $SL_{m,r}$ for each subinterval and then add the cumulative adjustments. 
    Since latter are usually relatively much smaller than the quadrature values, this might result in neglecting 
    the actual adjustment values. Please note that \texttt{Decimal} library was also used to address such constraints.
    \item [(11)] Ultimately, we add Monte Carlo result to the previous approximation.
\end{enumerate}

As it can be observed, the current solution enables the user to insert own interpolation meshes, increase 
the precision of computations, as well as extend the method to arbitrary regularity $r \in \mathbb{N}.$

\end{document}